# A note on the generalized-hypergeometric solutions of general and single-confluent Heun equations

D.Yu. Melikdzhanian[1,2] and A.M. Ishkhanyan[1,3]

[1]Russian-Armenian University, Yerevan, 0051 Armenia
[2]Vekua Institute of Applied Mathematics, Tbilisi State University, Tbilisi, 0186 Georgia
[3]Institute for Physical Research, Ashtarak, 0203 Armenia

We review the series solutions of the general and single-confluent Heun equations in terms of powers, ordinary-hypergeometric and confluent-hypergeometric functions. The conditions under which the expansions reduce to finite sums as well as the cases when the coefficients of power-series expansions obey two-term recurrence relations are discussed. Infinitely many cases for which a solution of a general or single-confluent Heun equation is written through a single generalized hypergeometric function are indicated. It is shown that the parameters of this generalized hypergeometric function are the roots of a certain polynomial equation, which we explicitly present for any given order.

## 1. Introduction

The five special functions of the Heun class [1-3] originating from the Heun class of ordinary differential equations suggest advanced mathematical objects that are currently widely encountered in physics and mathematics research. The applications include, e.g., classical and non-classical physics, mathematical physics, engineering, chemistry, statistics, mathematical economics, etc. [4]. Despite their importance, the theory of these functions at the present time needs essential development because these equations suggest several additional challenges not encountered in the case of their prominent predecessors, hypergeometric functions (for details, see [1] and references therein). The complications are caused by the extra singularity, as compared with the corresponding hypergeometric equations, which possess utmost three singular points [1-3].

In the vicinities of the regular singular points, the Heun equations admit solutions in the form of convergent power series – those by Frobenius [5]. In contrast, the power-series in the vicinities of irregular singular points are generally not convergent. Therefore, for the latter series of particular interest, first of all, are the cases when the series terminate thus turning into finite sums [6-9]. Besides, in several cases the Heun functions allow representations in the form of series expansions in terms of the hypergeometric functions of various types (see [10-23]). These series may or may not be convergent. For such series, the finite-sum reductions also play an important role because the properties of the hypergeometric functions are much more studied as compared with the Heun functions (see, e.g., [3,24,25]).

In the present paper, we discuss the finite-sum reductions of the hypergeometric-series expansions of the general and the single-confluent Heun functions. It is known that the result of such a reduction can be presented in terms of a generalized hypergeometric function [26-33]. A main result we report here is that the parameters of this generalized hypergeometric function are the roots of a certain polynomial equation of a corresponding degree. We present the explicit construction of this equation for the general case and write down several particular examples of it for the lowest order solutions.

## 2. General approach

We discuss the solution of a Heun second-order linear differential equation:

$$\frac{d^2\Phi}{dz^2} + G(z)\frac{d\Phi}{dz} + H(z)\Phi = 0, \qquad (1)$$

as a series in terms of simpler mathematical functions $u_n(z)$:

$$\Phi(z) = \sum c_n u_n(z) \qquad (2)$$

(in particular, this can be a power series for which $u_n = z^n$). In general, this series may be single- or double-sided, convergent, or asymptotic. The convergence region of the series depends on the particular equation under consideration and the type of the applied expansion functions. For instance, if power-series expansions are discussed, the convergence region for the general Heun equation is the disk of radius one while for the confluent Heun equation the expansion is convergent everywhere on the complex $z$-plane. The focus of our research is on the particular cases when the series terminates into a finite sum. The motivation for this is that the convergence region for such finite-sum solutions can be extended to the whole complex plane.

Furthermore, we restrict the discussion by the cases when the series coefficients $c_n$ obey a *three-term* recurrence relation between successive coefficients (note that there exist many expansions for which the coefficients satisfy four- or more-term recurrence relations, see, e.g., [21-23]). The reason is that the three-term relations are the most often encountered ones. Indeed, the power-series expansions [1-5] as well as the classical hypergeometric-series expansions [10-20] all obey three-term recurrence relations. Thus, we suppose that the series coefficients satisfy the relation

$$R_n c_n + Q_{n-1} c_{n-1} + P_{n-2} c_{n-2} = 0, \qquad (3)$$

where we assume $c_0 = 1$, without loss of the generality. In order the series be single-sided with the summation index $n$ running from zero to plus infinity, we assume $c_{-1} = c_{-2} = 0$ and

$$R_0 = 0, \quad R_n \neq 0, \quad n = 1, 2, 3... \qquad (4)$$

The series then terminates on the $(N+1)$th term, thus turning into the finite sum

$$\Phi(z) = \sum_{n=0}^{N} c_n u_n(z), \qquad (5)$$

if $c_N \neq 0$, $c_{N+1} = c_{N+2} = 0$. Equation (3) shows that these conditions are satisfied if

$$R_{N+1} \cdot 0 + Q_N c_N + P_{N-1} c_{N-1} = 0 \qquad (6)$$

and

$$R_{N+2} \cdot 0 + Q_{N+1} \cdot 0 + P_N c_N = 0. \qquad (7)$$

As a result, one arrives at the following conditions:

$$P_N = 0 \qquad (8)$$

and

$$\det \mathbf{M}^{(N)} = 0, \qquad (9)$$

where $\mathbf{M}^{(N)}$ is the $(N+1) \times (N+1)$ square minor formed by the first $(N+1)$ rows and the first $(N+1)$ columns of the infinite matrix

$$\mathbf{M} = \begin{bmatrix} Q_0 & R_1 & 0 & 0 & 0 & ... \\ P_0 & Q_1 & R_2 & 0 & 0 & ... \\ 0 & P_1 & Q_2 & R_3 & 0 & ... \\ 0 & 0 & P_2 & Q_3 & R_4 & ... \\ ... & ... & ... & ... & ... & ... \end{bmatrix}. \qquad (10)$$

It can be shown that in general equation (9) presents a $(N+1)$th degree polynomial equation for the *accessory* parameter of the Heun equation. We recall that unlike the other parameters involved in the equation, the accessory parameter is the only one which does not characterize the behavior of a solution in the vicinity of a singular point (see [1,2]).

**Expansions in terms of generalized hypergeometric functions.** It is known that the solutions of the Heun equations of certain types can be expanded as (finite or infinite, convergent or divergent) series in terms of generalized hypergeometric functions $\tilde{u}_n(z)$ and $\tilde{\tilde{u}}_n(z)$ given as

$$\tilde{u}_n(z) = {}_pF_q\left(\alpha_1+n,\alpha_2,...,\alpha_p;\gamma_1,...,\gamma_q;\omega z\right), \tag{11}$$

$$\tilde{\tilde{u}}_n(z) = {}_pF_q\left(\alpha_1,...,\alpha_p;\gamma_1-n,\gamma_2...,\gamma_q;\omega z\right) \tag{12}$$

(see, e.g., [26]-[34]). The non-negative integer numbers $p$ and $q$ as well as the values of parameters $\alpha_j$ and $\gamma_k$ depend on the particular type of the Heun equation under consideration. Here we apply a different set of expansion functions $u_n(z)$ given as

$$u_n(z) = z^n \frac{d^n}{dz^n} {}_pF_q\left(\alpha_1,...,\alpha_p;\gamma_1,...,\gamma_q;\omega z\right) \tag{13}$$

A useful observation for what follows is that there exist certain relations between functions $\tilde{u}_n(z), \tilde{\tilde{u}}_n(z)$ and $u_n(z)$. Indeed, using the differentiation formulas for generalized hypergeometric functions [25]

$$\frac{d^n}{dz^n}\left(z^{\alpha_1+n-1}u_0(\omega z)\right) = (\alpha_1)_n z^{\alpha_1-1} {}_pF_q\left(\alpha_1+n,\alpha_2,...,\alpha_p;\gamma_1,...,\gamma_q;\omega z\right), \tag{14}$$

$$\frac{d^n}{dz^n}\left(z^{\gamma_1-1}u_0(\omega z)\right) = (\gamma_1-n)_n z^{\gamma_1-n-1} {}_pF_q\left(\alpha_1,...,\alpha_p;\gamma_1-n,\gamma_2,...,\gamma_q;\omega z\right), \tag{15}$$

where $(...)_n$ is the Pochhammer symbol [3] and

$$u_0(z) = {}_pF_q\left(\alpha_1,...,\alpha_p;\gamma_1,...,\gamma_q;z\right), \tag{16}$$

it can be readily shown that

$$\tilde{u}_n(z) = \sum_{k=0}^{n} \frac{C_n^k}{(\alpha_1)_k} z^k \frac{d^k \tilde{u}_0(z)}{dz^k} = \sum_{k=0}^{n} \frac{C_n^k}{(\alpha_1)_k} u_k(z), \tag{17}$$

$$\tilde{\tilde{u}}_n(z) = \sum_{k=0}^{n} \frac{C_n^k}{(\gamma_1-n)_k} z^k \frac{d^k \tilde{\tilde{u}}_0(z)}{dz^k} = \sum_{k=0}^{n} \frac{C_n^k}{(\gamma_1-n)_k} u_k(z), \tag{18}$$

where $C_n^k = \binom{n}{k}$ are the binomial coefficients. Thus, each of functions $\tilde{u}_n(z)$ and $\tilde{\tilde{u}}_n(z)$ can be presented as a linear combination of the functions $u_0(z)$, $u_1(z)$, ..., $u_n(z)$.

Furthermore, any finite sum

$$\sum_{n=0}^{N} d_n u_n(z), \tag{19}$$

that is a linear combination of functions $u_0(z)$, $u_1(z)$, ..., $u_n(z)$, can be presented as a linear combination of functions $\tilde{u}_0(z)$, $\tilde{u}_1(z)$, ..., $\tilde{u}_n(z)$ or $\tilde{\tilde{u}}_0(z)$, $\tilde{\tilde{u}}_1(z)$, ..., $\tilde{\tilde{u}}_n(z)$. The calculations show that if

$$\sum_{n=0}^{\infty} d_n u_n(z) = \sum_{n=0}^{\infty} \tilde{d}_n \tilde{u}_n(z) = \sum_{n=0}^{\infty} \tilde{\tilde{d}}_n \tilde{\tilde{u}}_n(z), \qquad (20)$$

then
$$d_n = \sum_{j=0}^{\infty} \frac{C_{n+j}^n}{(\alpha_1)_n} \tilde{d}_{n+j} = \sum_{j=0}^{\infty} \frac{C_{n+j}^n}{(\gamma_1 - n - j)_n} \tilde{\tilde{d}}_{n+j}. \qquad (21)$$

**Solutions of a Heun equation in terms of a single generalized hypergeometric function.** It is known that the general and single-confluent Heun equations admit solutions in terms of a single generalized hypergeometric function [25-33]. These solutions can be derived in the following way. Using the differentiation formula for hypergeometric functions [25]

$$\frac{d}{dz}\left(z^\lambda u_0(z)\right) = \lambda z^{\lambda-1} {}_{p+1}F_{q+1}\left(\alpha_1, ..., \alpha_p, \lambda+1; \gamma_1, ..., \gamma_q, \lambda; \omega z\right), \qquad (22)$$

we have the relation

$$\lambda \cdot {}_{p+1}F_{q+1}\left(\alpha_1, ..., \alpha_p, \lambda+1; \gamma_1, ..., \gamma_q, \lambda; \omega z\right) = z^{-\lambda} z \frac{d}{dz}\left(z^\lambda u_0(z)\right) = \left(z\frac{d}{dz} + \lambda\right) u_0(z). \qquad (23)$$

Applying this formula for several times, we obtain

$$\lambda_1 \cdot ... \cdot \lambda_N \cdot {}_{p+N}F_{q+N}\left(\alpha_1, ..., \alpha_p, \lambda_1+1, ..., \lambda_N+1; \gamma_1, ..., \gamma_q, \lambda_1, ..., \lambda_N; \omega z\right)$$
$$= \left(z\frac{d}{dz} + \lambda_1\right)...\left(z\frac{d}{dz} + \lambda_N\right) u_0(z). \qquad (24)$$

Furthermore, it is possible to show (e.g., using the method of mathematical induction) that

$$z^n \frac{d^n}{dz^n} = \left(z\frac{d}{dz}\right)\left(z\frac{d}{dz} - 1\right)...\left(z\frac{d}{dz} - n + 1\right). \qquad (25)$$

With this, we conclude that if a solution of equation (1) is given as the finite sum

$$\Phi(z) = \sum_{n=0}^{N} d_n z^n \frac{d^n u_0(z)}{dz^n}, \qquad (26)$$

then this solution can be written as

$$\Phi(z) = \sum_{n=0}^{N} d_n \left(z\frac{d}{dz}\right)\left(z\frac{d}{dz} - 1\right)...\left(z\frac{d}{dz} - n + 1\right) u_0(z)$$
$$= d_N \left(z\frac{d}{dz} + e_1\right)...\left(z\frac{d}{dz} + e_N\right) u_0(z), \qquad (27)$$

where $e_1, ..., e_N$ are the zeros, multiplied by $(-1)$, of the polynomial of degree $N$

$$A(\xi) = \sum_{n=0}^{N} d_n \xi(\xi-1)...(\xi-n+1). \qquad (28)$$

Or, otherwise, $e_1, ..., e_N$ are the roots of the algebraic equation

$$\sum_{n=0}^{N}(-1)^n d_n e(e+1)...(e+n-1) = 0. \qquad (29)$$

With this, taking into the account the equation (24) and omitting the insignificant constant factor, we arrive at the solution

$$\Phi(z) = {}_{N+p}F_{N+q}\left(\alpha_1,...,\alpha_p,e_1+1,...,e_N+1;\gamma_1,...,\gamma_q,e_1,...,e_N;\omega z\right). \qquad (30)$$

Equation (29) defining the parameters of this solution is a main result of the present paper. A complementary remark is that the coefficients of polynomial $A(\xi)$ are explicitly presented in terms of the Stirling numbers of the first kind $s_n^k$ [3]. Indeed,

$$A(\xi) \equiv d_0 + \sum_{k=1}^{N} A_k \xi^k = d_0 + \sum_{n=1}^{N} d_n \sum_{k=1}^{n}(-1)^{n-k} s_n^k \xi^k, \qquad (31)$$

so that

$$A_k = \sum_{n=k}^{N}(-1)^{n-k} s_n^k d_n. \qquad (32)$$

We conclude this section by noting that if a solution of a Heun equation presents a linear combination of a finite number of generalized hypergeometric functions, then the coefficients of the *Frobenius expansion* of this solution can conveniently be written in terms of polynomial $A(\xi)$. Indeed, since the solution $\Phi(z)$ can be presented in terms of a *single* generalized hypergeometric function (that is, equation (30)), hence, the coefficients of the power-series expansion of this function are given by the simple explicit expression

$$c_n = \frac{\omega^n (\alpha_1)_n \cdot ... \cdot (\alpha_p)_n}{n!(\gamma_1)_n \cdot ... \cdot (\gamma_q)_n} \cdot \frac{(e_1+n) \cdot ... \cdot (e_n+n)}{e_1 \cdot ... \cdot e_n} = \frac{\omega^n (\alpha_1)_n \cdot ... \cdot (\alpha_p)_n}{n!(\gamma_1)_n \cdot ... \cdot (\gamma_q)_n} \cdot \frac{A(n)}{A(0)}. \qquad (33)$$

Accordingly, in this case the coefficients $c_n$ satisfy the two-term recurrence relation

$$\frac{c_{n+1}}{c_n} = \frac{\omega(\alpha_1+n) \cdot ... \cdot (\alpha_p+n)}{(n+1)(\gamma_1+n) \cdot ... \cdot (\gamma_q+n)} \cdot \frac{A(n+1)}{A(n)}. \qquad (34)$$

Since the coefficients $c_n$ of a Frobenius series solution of the Heun equation (1) generally satisfy a three-term recurrence relation [1-3], the general conclusion is that if a solution of a Heun equation presents a linear combination of a finite number of generalized hypergeometric functions, then the three-term recurrence relation for the coefficients of the Frobenius expansion of this solution reduces to a two-term one. This observation served as a starting point for construction of the generalized hypergeometric solutions of the general and single-confluent Heun equations in [**Error! Reference source not found.**,33].

## 3. An expansion of the general Heun function in terms of the Gauss functions

Consider the general Heun equation

$$\frac{d^n}{dz^n}\Phi(z)+\left(\frac{\gamma}{z}+\frac{\delta}{z-1}+\frac{\varepsilon}{z-a}\right)\frac{d}{dz}\Phi(z)+\frac{\alpha\beta z-q}{z(z-1)(z-a)}\Phi(z)=0, \quad (35)$$

where the parameters obey the Fuchsian condition

$$\gamma+\delta+\varepsilon=\alpha+\beta+1. \quad (36)$$

Let the solution of equation (35) is expanded in terms of the following functions involving the Gauss ordinary hypergeometric functions:

$$u_n(z)=z^n\frac{d^n}{dz^n}{}_2F_1(\alpha,\beta;\gamma;z), \quad (37)$$

which differ from the ones applied in [14]. These functions satisfy the differential equation

$$\frac{d^n u_n}{dz^n}+\left(\frac{\gamma-n}{z}+\frac{\delta+n}{z-1}\right)\frac{du_n}{dz}+\left(\frac{n(1-\gamma)}{z^2}+\frac{\alpha\beta+n(\alpha+\beta-\delta)}{z(z-1)}\right)u_n=0, \quad (38)$$

and the three-term recurrence relation

$$(z-1)u_n(z)+\left((\alpha+\beta+2n-3)z-\gamma-n+2\right)u_{n-1}(z)+z(n+\alpha-2)(n+\beta-2)u_{n-2}(z)=0. \quad (39)$$

With this, it is checked that the coefficients $d_n$ of the expansion

$$\Phi(z)=\sum_{n=0}^{\infty}d_n u_n(z) \quad (40)$$

satisfy the three-term recurrence relation

$$R_n d_n+Q_{n-1}d_{n-1}+P_{n-2}d_{n-2}=0 \quad (41)$$

with coefficients $R_n, Q_n, P_n$ given as

$$R_n=an(n+\alpha-1)(n+\beta-1), \quad (42)$$

$$Q_n=a(n+\alpha)(n+\beta)+(a-1)n(n+\varepsilon-1)-\gamma n-q, \quad (43)$$

$$P_n=(a-1)(n+\varepsilon). \quad (44)$$

Consider the cases when the series terminates into a finite sum in which the maximum number of non-zero coefficients is equal to $N+1$ (that is, the last non-zero coefficient is $d_N$). For this to happen, the necessary and sufficient conditions are

$$\varepsilon=-N \quad (45)$$

and $q$ is an eigenvalue of the $(N+1)\times(N+1)$ square minor of the matrix

$$\begin{bmatrix} \tilde{Q}_0 & a\alpha\beta & 0 & 0 & ... \\ (a-1)\varepsilon & \tilde{Q}_1 & 2a(\alpha+1)(\beta+1) & 0 & ... \\ 0 & (a-1)(\varepsilon+1) & \tilde{Q}_2 & 3a(\alpha+2)(\beta+2) & ... \\ 0 & 0 & (a-1)(\varepsilon+2) & \tilde{Q}_3 & ... \\ ... & ... & ... & 0 & ... \end{bmatrix}, \quad (46)$$

where $\tilde{Q}_n = Q_n + q$. Note that the coefficients $d_0$, ..., $d_N$ are then the components of corresponding eigenvector.

In particular, at $N = 0$ the conditions explicitly read

$$\varepsilon = 0, \quad q = a\alpha\beta, \quad (47)$$

at $N = 1$ we have

$$\varepsilon = -1, \quad (48)$$

$$q^2 - (1 - \gamma + a(2\alpha\beta + \alpha + \beta))q + a\alpha\beta(-\gamma + a(\alpha+1)(\beta+1)) = 0, \quad (49)$$

at $N = 2$ should be

$$\varepsilon = -2, \quad (50)$$

$$\varsigma^3 + (a(3\alpha + 3\beta + 1) + 4 - 3\gamma)\varsigma^2$$
$$+ (2(a(\alpha+\beta-1) + 2 - \gamma)(a(\alpha+\beta+1) + 1 - \gamma) + 2a(a-1)(2\alpha\beta + \alpha + \beta + 1))\varsigma, \quad (51)$$
$$+ 4(a-1)a\alpha\beta((a(\alpha+\beta+1) + 1 - \gamma)) = 0$$

where

$$\varsigma = a\alpha\beta - q. \quad (52)$$

We note that, according to the discussion in Section 2, the expansion of a solution $\Phi(z)$ of the general Heun equation (35) in terms of functions $u_n(z)$ given by equation (37) can readily be transformed to expansion of this solution in terms of functions

$$\tilde{u}_n(z) = {}_2F_1(\alpha + n, \beta; \gamma, z) \quad (53)$$

or

$$\tilde{\tilde{u}}_n(z) = {}_2F_1(\alpha, \beta; \gamma - n, z). \quad (54)$$

In other words, if the solution $\Phi(z)$ can be presented as a linear combination of a finite number of functions $u_n(z)$, then this solution can be presented as a linear combination of functions $\tilde{u}_n(z)$ or $\tilde{\tilde{u}}_n(z)$ as well. The expansion of the Heun function in terms of the Gauss hypergeometric functions $\tilde{\tilde{u}}_n(z)$, for which the summation index is involved only in the denominator parameter of the Gauss functions, has been discussed in [14], where one can find the corresponding three-term recurrence relation for coefficients of that expansion.

**Solutions of the general Heun equation in terms of generalized hypergeometric functions.** According to the general discussion presented in Section 2, if a solution of equation (35) can be presented as the linear combination

$$\Phi(z) = \sum_{n=0}^{N} d_n u_n(z) = \sum_{n=0}^{N} d_n z^n \frac{d^n}{dz^n} {}_2F_1(\alpha,\beta;\gamma;z), \quad (55)$$

this solution can also be written as the generalized hypergeometric function

$$\Phi(z) = {}_{2+N}F_{1+N}\left(\alpha,\beta,e_1+1,...,e_N+1;\gamma,e_1,...,e_N;z\right), \quad (56)$$

where $e_1, ..., e_N$ are the roots of polynomial equation (29). Here are the explicit solutions of the three lowest orders $N = 0,1,2$ (note that these solutions apply for $\varepsilon, q$ defined by equations (47) if $N = 0$, (48),(49) if $N = 1$, and (50)-(52) if $N = 2$).

$N = 0$: 
$$\Phi(z) = {}_2F_1(\alpha,\beta;\gamma;z). \quad (57)$$

$N = 1$: 
$$\Phi(z) = {}_3F_2\left(\alpha,\beta,e_1+1;\gamma,e_1;z\right), \quad (58)$$

with 
$$e_1 = \frac{1}{d_1} = -\frac{1}{1-q/(a\alpha\beta)}. \quad (59)$$

$N = 2$: 
$$\Phi(z) = {}_4F_3\left(\alpha,\beta,e_1+1,e_2+1;\gamma,e_1,e_2;z\right), \quad (60)$$

where $e_1$ and $e_2$ are the roots of the quadratic equation

$$d_2 e(e+1) - d_1 e + d_0 = 0 \quad (61)$$

with 
$$d_0 = 1, \quad d_1 = -1 + q/(a\alpha\beta), \quad (62)$$

$$d_2 = -\frac{d_1}{2} + \frac{2(a-1)(1+d_1)+d_1(\gamma+q)}{2a(\alpha+1)(\beta+1)} = \frac{1}{2} - \frac{q}{2a\alpha\beta} + \frac{2(a-1)q+(q-a\alpha\beta)(\gamma+q)}{2a^2\alpha\beta(\alpha+1)(\beta+1)}. \quad (63)$$

We conclude this section by noting that if the solution (56) is expanded as a power series, then the coefficients of this expansion are explicitly written as

$$c_n = \frac{(\alpha)_n(\beta)_n}{n!(\gamma)_n} \cdot \frac{(e_1+n)\cdot...\cdot(e_n+n)}{e_1\cdot...\cdot e_n} = \frac{(\alpha)_n(\beta)_n}{n!(\gamma)_n} \cdot \frac{A(n)}{A(0)} \quad (64)$$

where $A(\xi)$ is the polynomial defined by equation (28) or (31). Note that these coefficients satisfy the two-term recurrence relation

$$\frac{c_{n+1}}{c_n} = \frac{(\alpha+n)(\beta+n)}{(n+1)(\gamma+n)} \cdot \frac{A(n+1)}{A(n)}. \quad (65)$$

These results have been derived in [**Error! Reference source not found.**] using a different technique. The advantage of the method applied here is that it directly leads to the explicit construction of the polynomial equation the zeros of which present the parameters of the generalized-hypergeometric solution of the general Heun equation. Besides, as already mentioned above, it is a systematic technique that can be applied to treat other Heun equations as well.

### 3. An expansion of the single-confluent Heun function in terms of the Kummer functions

The last statement can be supported by discussing the expansion of the solution of the single-confluent Heun equation

$$\frac{d^2}{dz^2}\Phi(z)+\left(\varepsilon+\frac{\gamma}{z}+\frac{\delta}{z-1}\right)\frac{d}{dz}\Phi(z)+\frac{\varepsilon\alpha z-q}{z(z-1)}\Phi(z)=0 \qquad (66)$$

in terms of the following functions involving the Kummer confluent hypergeometric functions:

$$u_n(z)=z^n\frac{d^n}{dz^n}{}_1F_1(\alpha,\gamma,-\varepsilon z). \qquad (67)$$

Note that the coefficient in front of $z$ in the numerator of the last term in equation (66) is put $\varepsilon\alpha$, not just $\alpha$ as it is put in some references. Note also that because of this convention the expansion does not apply to the Ince limit $\varepsilon=0$ [35]. For discussion of the latter case the reader is referred to [33].

The expansion functions (67) differ from the ones applied in [15]. To construct an expansion of the single-confluent Heun function in terms of functions (67), we note that they satisfy the differential equation

$$\frac{d^n u_n}{dz^n}+\left(\frac{\gamma-n}{z}+\varepsilon\right)\frac{du_n}{dz}+\frac{n(1-\gamma)+\varepsilon\alpha z}{z^2}u_n=0, \qquad (68)$$

and three-term recurrence relation

$$\varepsilon z(\alpha+n-1)u_{n-1}(z)+(\gamma+n-1+\varepsilon z)u_n(z)+u_{n+1}(z)=0. \qquad (69)$$

With this, it is shown that the coefficients $d_n$ of the expansion

$$\Phi(z)=\sum_{n=0}^{\infty}d_n u_n(z) \qquad (70)$$

obey the three-term recurrence relation

$$R_n d_n+Q_{n-1}d_{n-1}+P_{n-2}d_{n-2}=0, \qquad (71)$$

the coefficients of which are given as

$$R_n=\varepsilon n(n+\alpha-1), \qquad (72)$$

$$Q_n=n(n+\varepsilon+\gamma+\delta-1)+\varepsilon\alpha-q, \qquad (73)$$

$$P_n = n + \delta. \tag{74}$$

The series (70) turns into a finite sum with the last non-zero coefficient being $d_N$ if

$$\delta = -N \tag{75}$$

and $(q - \alpha\varepsilon)$ is an eigenvalue of the $(N+1)\times(N+1)$ order minor of the matrix

$$\begin{bmatrix} 0 & \varepsilon\alpha & 0 & 0 & \dots \\ \delta & \varepsilon+\gamma+\delta & 2\varepsilon(\alpha+1) & 0 & \dots \\ 0 & \delta+1 & 2(\varepsilon+\gamma+\delta+1) & 3\varepsilon(\alpha+2) & \dots \\ 0 & 0 & \delta+2 & 3(\varepsilon+\gamma+\delta+2) & \dots \\ \dots & \dots & \dots & 0 & \dots \end{bmatrix}. \tag{76}$$

Note that the coefficients $d_0$, ..., $d_N$ are then the components of the corresponding eigenvector.

In particular, at $N = 0$ it should be

$$\delta = 0, \quad q - \varepsilon\alpha = 0. \tag{77}$$

Similarly, at $N = 1$ it should hold

$$\delta = -1, \tag{78}$$

$$(q - \varepsilon\alpha)^2 - (\varepsilon+\gamma-1)(q-\varepsilon\alpha) + \varepsilon\alpha = 0, \tag{79}$$

and at $N = 2$

$$\delta = -2, \tag{80}$$

$$(q-\varepsilon\alpha)^3 - (3\varepsilon+3\gamma-4)(q-\varepsilon\alpha)^2 + \\ \left(2(\varepsilon+\gamma-1)(\varepsilon+\gamma-2) + 2\varepsilon(2\alpha+1)\right)(q-\varepsilon\alpha) - 4\alpha\varepsilon(\varepsilon+\gamma-1) = 0. \tag{81}$$

We note that, as in the previous case of the general Heun equation, the expansion of a solution $\Phi(z)$ of the single-confluent Heun equation (66) in terms of functions $u_n(z)$ given by equation (67) can readily be transformed to expansion of this solution in terms of functions

$$\tilde{u}_n(z) = {}_1F_1(\alpha+n;\gamma;-\varepsilon z) \tag{82}$$

or

$$\tilde{\tilde{u}}_n(z) = {}_1F_1(\alpha;\gamma-n;-\varepsilon z). \tag{83}$$

Thus, if the solution $\Phi(z)$ can be presented in the form of a finite linear combination of the functions $u_n(z)$, then this solution can be presented in the form of linear combinations of functions $\tilde{u}_n(z)$ or $\tilde{\tilde{u}}_n(z)$. The expansion of the solution $\Phi(z)$ in terms of the functions $\tilde{u}_n(z)$ is presented in [15].

**Solutions of the single-confluent Heun equation in the terms of the generalized**

**confluent hypergeometric functions.** According to the general discussion presented in Section 2, the finite-term solution of equation (66)

$$\Phi(z) = \sum_{n=0}^{N} d_n u_n(z) = \sum_{n=0}^{N} d_n z^n \frac{d^n}{dz^n} {}_1F_1(\alpha; \gamma; -\varepsilon z) \qquad (84)$$

can be presented through a single generalized hypergeometric function written as

$$\Phi(z) = {}_{N+1}F_{N+1}(\alpha, e_1+1, ..., e_N+1; \gamma, e_1, ..., e_N; -\varepsilon z), \qquad (85)$$

where $e_1, ..., e_N$ are the roots of equation (29). We note that then the coefficients of the power-series expansion of this solution are given as

$$c_n = (-\varepsilon)^n \frac{(\alpha)_n}{n!(\gamma)_n} \cdot \frac{(e_1+n)\cdot...\cdot(e_n+n)}{e_1\cdot...\cdot e_n} = (-\varepsilon)^n \frac{(\alpha)_n(\beta)_n}{n!(\gamma)_n} \cdot \frac{A(n)}{A(0)}, \qquad (86)$$

where $A(\xi)$ is the polynomial defined by equation (28) or (31). Accordingly, these coefficients satisfy the two-term recurrence relation

$$\frac{c_{n+1}}{c_n} = -\varepsilon \frac{(\alpha+n)}{(n+1)(\gamma+n)} \cdot \frac{A(n+1)}{A(n)}. \qquad (87)$$

The first three explicit solutions read (note that these solutions apply for $\delta, q$ defined by equations (77) if $N=0$, (78),(79) if $N=1$, and (80),(81) if $N=2$).

$N=0$: $\qquad \Phi(z) = {}_1F_1(\alpha; \gamma; -\varepsilon z), \qquad (88)$

$N=1$: $\qquad \Phi(z) = {}_2F_2(\alpha, e_1+1; \gamma, e_1; -\varepsilon z), \qquad (89)$

where $\qquad e_1 = \dfrac{1}{d_1} = \dfrac{1}{-1+q/(\varepsilon\alpha)}, \qquad (90)$

$N=2$: $\qquad \Phi(z) = {}_3F_3(\alpha, e_1+1, e_2+1; \gamma, e_1, e_2; -\varepsilon z), \qquad (91)$

where $e_1$ and $e_2$ are the roots of the quadratic equation

$$d_2 e(e+1) - d_1 e + d_0 = 0 \qquad (92)$$

with $\qquad d_0 = 1, \quad d_1 = -\dfrac{(\varepsilon\alpha - q)d_0}{\varepsilon\alpha} = -1 + \dfrac{q}{\varepsilon\alpha}, \qquad (93)$

$$d_2 = -\frac{(\gamma+\delta+\varepsilon+\varepsilon\alpha-q)d_1+\delta d_0}{2\varepsilon(\alpha+1)} = \frac{q^2 - q(\gamma+\delta+\varepsilon+2\varepsilon\alpha) + \varepsilon\alpha(\gamma+\varepsilon+\varepsilon\alpha)}{2\varepsilon^2\alpha(\alpha+1)}. \qquad (94)$$

The generalized hypergeometric solution (85) of the Heun single-confluent equation has been discussed in [33]. As already mentioned above, the new point we report here is the polynomial equation (29) for parameters $e_1,...,e_N$. To give one more example of application of this equation, here is the solution of the single-confluent Heun equation for $\delta = -3$ ($N = 3$)

$$\Phi(z) = {}_4F_4(\alpha, e_1+1, e_2+1, e_3+1; \gamma, e_1, e_2, e_3; -\varepsilon z), \tag{95}$$

where $e_1, e_2, e_3$ are the roots of the cubic equation

$$-d_3 e(e+1)(e+2) + d_2 e(e+1) - d_1 e + d_0 = 0, \tag{96}$$

the coefficients $d_{0,1,2,3}$ being defined by the recurrence (71)-(74). This solution applies if the accessory parameter $q$ satisfies the fourth-order polynomial equation

$$d_4 = 0. \tag{97}$$

**Acknowledgments**

This research was supported by the Russian-Armenian (Slavonic) University at the expense of the Ministry of Education and Science of the Russian Federation, the Armenian Science Committee (SC Grant 20RF-171), and the Armenian National Science and Education Fund (ANSEF Grant No. PS-5701).